\documentclass[12pt,a4paper]{amsart}
\usepackage{amsfonts,amsmath,amssymb}

\def\SU{\operatorname{\mathrm {SU}}\nolimits}
\def\su{\operatorname{\mathfrak {su}}\nolimits}

\def\U{\operatorname{\mathrm {U}}\nolimits}
\def\SL{\operatorname{\mathrm {SL}}\nolimits}
\def\sl{\operatorname{\mathfrak {sl}}\nolimits}
\def\Sp{\operatorname{\mathrm {Sp}}\nolimits}
\def\sp{\operatorname{\mathfrak {sp}}\nolimits}
\def\Lie{\operatorname{\mathrm {Lie}}\nolimits}
\def\Aut{\operatorname{\mathrm {Aut}}\nolimits}

\def\diag{\operatorname{\mathrm {diag}}\nolimits}
\def\vol{\operatorname{\mathrm {Vol}}\nolimits}
\def\sign{\operatorname{\mathrm {sign}}\nolimits}
\def\Heis{\operatorname{\mathrm {Heis}}\nolimits}
\def\Ind{\operatorname{\mathrm {Ind}}\nolimits}
\def\ind{\operatorname{\mathrm {Ind}}\nolimits}

\def\cohind{\operatorname{\mathrm {Coh-Ind}}\nolimits}
\def\tr{\operatorname{\mathrm {tr}}\nolimits}
\def\ad{\operatorname{\mathrm {ad}}\nolimits}
\def\Ad{\operatorname{\mathrm {Ad}}\nolimits}
\def\sgn{\operatorname{\mathrm {sign}}\nolimits}
\def\Cent{\operatorname{\mathrm {Cent}}\nolimits}
\def\Gal{\operatorname{\mathrm {Gal}}\nolimits}
\def\bbR{\mathbb R}
\def\bbC{\mathbb C}
\def\bbD{\mathbb D}
\def\bbH{\mathbb H}
\def\bbN{\mathbb N}
\def\bbP{\mathbb P}
\def\bbS{\mathbb S}
\def\bbT{\mathbb T}

\def\beginpf{\textsc{Proof.}}
\def\endpf{\hfill$\Box$}
\newtheorem{thm}{Theorem}[section]
\newtheorem{prop}[thm]{Proposition}

\title[Poisson Summation and Endoscopy for $\Sp(4,\mathbb R)$]{Poisson Summation and Endoscopy\\ for $\Sp(4,\mathbb R)$}
\author[Do Ngoc Diep, Do Thi Phuong Quynh]{Do Ngoc Diep${}^{1}$, Do Thi Phuong Quynh${}^{2}$}
\address{${}^{1}$ Institute of Mathematics, VAST, 18 Hoang Quoc Viet road, Cau Giay district, 10307 Hanoi, Vietnam \newline
{\tt Email: dndiep@math.ac.vn}}
\address{${}^{2}$ {\sc Medicine University Colllege, Thai Nguyen Universtiy, Thai Nguyen City, Vietnam}\newline
{\tt Email:  phuongquynhtn@gmail.com}}
\begin{document}
\date{\bf Version of \today}
%
\maketitle
\begin{abstract} 
In this paper we analyze the endoscopy for $\Sp(4,\mathbb R)$. The new results are a precise realization of the discrete series representations (in Section 2), a computation of their traces (Section 3) and an exact formula for the noncommutative Poisson summation and endoscopy of  for this group (in Sections 4,5).
\\
\textit{Key terms}: trace formula; orbital integral; transfer; endoscopy
\end{abstract}
\tableofcontents
\section{Introduction}
The following arguments were introduced in the previous work \cite{diepquynh2}: 
For a fixed reductive group  we do find, if possible all the irreducible unitary representations  and then decompose  a particular represention into a discrete or continuous direct sum or integral of irreducible ones  are the basic questions of harmonic analysis on reductive Lie groups. In particlular, the discrete part of the regular representation of reductive groups is the discrete sum of discrete series representations. 

Often these problems are reduced to the trace formula, because, as known, the unitary representations are uniquely defined by it generalized character and infinitesimal character. Following Harish-Chandra, the generalized character is defined by its restriction to the maximal compact subgroup, as the initial eigenvalue problem for the generalized Laplacian (the Casimir operator) with the infinitesimal character as the eigenvalues infinitesimal action of Casimir operators.
There is a very highly developed theory of Arthur-Selberg trace formula. The theory is complicated and one reduces it to the same problem for smaller   endoscopic subgroups. It is called the transfer and plays a very important role in the theory.
By definition, an endoscopic subgroup is the connected component of the centralizer of regular semisimple elements, associated to representations, namely by the orbit method.

In the previous papers \cite{diepquynh1}-\cite{diepquynh3}, we treated the case of rank one $\SL(2,\bbR)$ and  $\SL(3,\bbR)$, $\SU(2,1)$ of rank 2. In this paper we a doing the same study for $\Sp(4,\bbR)$ also of rank 2.

For the discrete series representations of $\Sp(4,\bbR)$ the following endoscopic groups should be considered:
\begin{itemize}
\item \textit{the elliptic case}: diagonal subgroup of regular elements
$$H= \left\{\diag(a_1,a_2,a_1^{-1}, a_2^{-1}); \quad a_1\ne a_2,\right\};$$
\item \textit{the parabolic case}: blog-diagonal subgroup of regular elements 
 $\gamma = k_{\theta_1}k_{\theta_2} = \begin{pmatrix}\cos\theta_1 & 0 &\sin\theta_1 & 0\\ 0 & \cos\theta_2 &0 & \sin\theta_2 \\ -\sin\theta_1 &0 &  \cos\theta_1 & 0\\  0 & -\sin\theta_2 & 0 &\cos\theta_2 \end{pmatrix}$. 
\item the trival case: the group $H= \Sp(4,\bbR)$ its-self.
\end{itemize}
We show in this paper that the method that J.-P. Labesse \cite{labesse} used for $\SL(2,\bbR)$ is applied also for $\Sp(4,\bbR)$.
Therefore one deduces the transfer formula for the discrete series representations and limits of  $\Sp(4,\bbR)$ to the corresponding endoscopic group. 

For the group $\Sp(4,\bbR)$ we make a precise realization of the discrete series representations (in Section 2) by using the Orbit Method and Geometric Quantization to the solvable radical, a computation in the context of $\Sp(4,\bbR)$ of their traces (Section 3) and an exact formula for the noncommutative Poisson summation and endoscopy of  for this group (in Sections 4,5). 

\section{Irreducible Unitary Representations of $\Sp(4,\mathbb R)$}

\subsection{The structure of $\Sp(4,\mathbb R)$} The following notions and results are folklore and we recall them to fix an appropriate system of notations.

Let us remind that the group $\Sp(4,\mathbb R)$ is 
$$\Sp(4,\mathbb R)= \left\{   \begin{pmatrix} A & B \\ C & D \end{pmatrix} \left| {A,B,C,D \in M_2(\mathbb R)\atop
\begin{pmatrix} A & B \\ C & D \end{pmatrix} J_4 + J_4 \begin{pmatrix} A & B \\ C & D \end{pmatrix} = 0 }\right. \right\}$$ where $J$ is the matrix of the skew symmetrix form
$$\omega(u,v) = \left| \begin{matrix} u_1 & u_2\\ v_1 & v_2  \end{matrix} \right| + \left| \begin{matrix} u_3 & u_4\\ v_3 & v_4  \end{matrix} \right| , u,v \in \mathbb R^4, $$ i.e. 
$$J_4 = \begin{pmatrix} 0 & 1 & 0 & 0\\ -1 &0  & 0 & 0 \\ 0 & 0&0 & 1\\ 0 & 0 & -1 & 0 \end{pmatrix}, \omega(u,v) = {}^tuJ_4v.$$

Denote by $\sp(4,\mathbb R)$ the Lie algebra $\Lie\Sp(4,\mathbb R)$, $\theta$ the Cartan involution of the group $G=\Sp(4,\mathbb R)$. 
The corresponding Cartan involution for its Lie algebra $\sp(4,\mathbb R)$  is denote by the same symbol $\theta\in \Aut\sp(4,\mathbb R) $,  $$\theta(X) = {}^tX^{-1}, X\in \sp(4,\mathbb R)$$ 
Therefore, the Lie algebra $\sp(4,\mathbb R)$ can be described as
$$\sp(4,\mathbb R)= \left\{ \begin{pmatrix} A & B \\ C & D \end{pmatrix} \left| {A,B,C,D \in M_2(\mathbb R)\atop {}^tA = -D, {}^tB = B, {}^tC=C }  \right.\right\}$$
The maximal compact subgroup $K$ of $G$
$$K = \left\{\left. \begin{pmatrix} A & B\\ -B & A\end{pmatrix} \right|A,B \in M_2(\mathbb R)    \right\} \cong U(2) = \left\{ U=A + iB | {}^t\bar{U}U=I_2\right\} $$ 
 is the subgroup of $G$, the Lie algebra $\mathfrak k$ of which is consisting of all the matrices with eigenvalue $+1$,
$$\mathfrak k =\Lie K  =\left\{\left. \begin{pmatrix} A & B\\ -B & A\end{pmatrix}  \right|  {}^tA = -A, {}^tB = B   \right\}.$$

The \textsl{Borel subgroup} of $\Sp(4,\mathbb R)$ is the \textsl{minimal parabolic} subgroup $$P_0 = B=AU, A =\left\{\diag(t_1,t_2,t_2^{-1},t_1^{-1}) | t_1,t_2 > 0 \right\},$$
$$U = \left\{ n(x_1,x_2,x_3,x_4) = 
\begin{pmatrix} 1 & 0 & x_1 & x_2\\ 0 & 1 & x_2 & x_3 \\
0 & 0 & 1 & 0\\
0 & 0 & 0 & 1
 \end{pmatrix}\begin{pmatrix}  
1 & x_4 & 0 & 0\\ 0 & 1 & 0 & 0\\ 0 & 0& 1 & 0\\ 0 & 0 & -x_4 & 1
\end{pmatrix} \right\},$$
the Lie algebra of which is consisting of all the matrices with eigenvalue $-1$,
$$\mathfrak b = \mathfrak u + \mathfrak a=\left\{\left.  X\in \mathfrak g \right| \theta(X) = -X   \right\}$$
$$\mathfrak p = \left\{  \begin{pmatrix} * & * & *&  *\\ 0 & * & * & * \\ 0 & 0 & * & 0\\ 0 & 0 & * & *  \end{pmatrix}\right\} = \langle H_1,H_2,E_{2\mathbf e_1},E_{\mathbf e_1 + \mathbf e_2}, E_{2\mathbf e_2}, E_{\mathbf e_1-\mathbf e_2} \rangle,$$  where
$$H_1 =  \begin{pmatrix} 1 & 0 & 0 & 0\\
0 & 0 & 0 & 0\\ 0 & 0 &-1 & 0\\ 0 & 0 & 0 & 0\end{pmatrix}, H_2  = \begin{pmatrix} 0 & 0 & 0 & 0\\
0 & 1 & 0 & 0\\ 0 & 0 &0 & 0\\ 0 & 0 & 0 & -1\end{pmatrix}$$
$$E_{2\mathbf e_1} = \begin{pmatrix} 0& 0 & 1 & 0\\  0 & 0 & 0 & 0\\ 0 & 0& 0& 0\\ 0 & 0& 0 & 0\end{pmatrix},\quad  E_{\mathbf e_1+\mathbf e_2} = \begin{pmatrix} 0& 0 & 0 & 1\\  0 & 0 & 1 & 0\\ 0 & 0& 0& 0\\ 0 & 0& 0 & 0\end{pmatrix}$$
$$E_{2\mathbf e_2} = \begin{pmatrix} 0& 0 & 0 & 0\\  0 & 0 & 0 & 1\\ 0 & 0& 0& 0\\ 0 & 0& 0 & 0\end{pmatrix},\quad  E_{\mathbf e_1-\mathbf e_2} = \begin{pmatrix} 0& 1 & 0 & 0\\  0 & 0 & 0 & 0\\ 0 & 0& 0& 0\\ 0 & 0& -1 & 0\end{pmatrix}$$

The other parabolic subgroup is the \textsl{Jacobi parabolic} subgroup $P_J= M_JA_JU_J$,
$A_J = \left\{\diag(1,t_1,t_1^{-1},1) | t_1 > 0 \right\}\cong \mathbb R_+^*,$
$U_J = \left\{ n(x_1,x_2,0,x_4) | x_1,x_2,x_4 \in \mathbb R \right\}\cong \Heis(3,\mathbb R),$
$M_J = \SL(2,\mathbb R) \times \{\pm1\},$
the  Lie algebra of which is $$\mathfrak p_J = \mathfrak u_J \oplus \mathfrak a_J \oplus \mathfrak m_J = \left\{ \begin{pmatrix} * & * & * & * \\ 0 & * & * & *\\ 0 & * & * & * \\ 0 & 0 & 0 & * \end{pmatrix}  \right\},$$
where $$\mathfrak a_J = \langle H_1 \rangle, \mathfrak m_J = \langle H_2, T_2,  E_{2\mathbf e_2} \rangle \cong \sl(2,\mathbb R), \mathfrak u_J = \langle E_{2\mathbf e_1}, E_{\mathbf e_1 + \mathbf e_2}, E_{\mathbf e_1 - \mathbf e_2} \rangle \cong \Heis(3)$$

The group $\Sp(4,\mathbb R)$ admits the well-known Cartan decomposition in form of a semi-direct product
$G= B\rtimes K$, denoted simply by $BK$. 

The complexified Lie algebra $$\mathfrak g_\mathbb C  = \sp(4,\mathbb C)= \mathfrak h \oplus \sum_{\beta\ne 0} \mathfrak g_\beta =\mathfrak h \oplus \mathfrak p_+ \oplus \mathfrak p_-$$  
$$\mathfrak p_+= \bigoplus_{\beta > 0} \mathfrak g_\beta,
\mathfrak p_-= \bigoplus_{\beta < 0} \mathfrak g_\beta,$$
where $\mathfrak g_\beta = \langle X_\beta \rangle$ is the root space corresponding to the root $\beta$.

The associate root system is 
$$\Sigma = \left\{\pm(2,0), \pm(0,2), \pm(1,1),\pm(1,-1)  \right\}$$ from which the compact positive  roots are $\Delta_c =\{(1,-1) \}$ and the noncompact positive roots are $\Delta_n =\{(2,0),(1,1),(0,2)\}$. The root vectors are 
$$X_{\pm(2,0)} = \begin{pmatrix} 
1 & 0 &\pm i & 0\\
0 & 0 & 0 & 0 \\
\pm i& 0 &-1& 0\\
0 & 0 & 0 & 0
\end{pmatrix}, X_{\pm(0,2)} = \begin{pmatrix} 
0 & 0 &0 & 0\\
0 & 1 & 0 & \pm i \\
0 & 0 & 0 & 0\\
0 & \pm i & 0 & 1
\end{pmatrix}$$ 
$$X_{\pm(1,1)} = \begin{pmatrix} 
0 & 1 &0  & \pm i\\
1 & 0 & \pm i & 0 \\
0 & \pm i &0 & -1\\
\pm i & 0 &-1 & 0
\end{pmatrix}, X_{\pm(1,-1)} = \begin{pmatrix} 
0 & 1 &0 &\pm i\\
-1 & 0  & \pm i & 0 \\
0 &  \pm i & 0 &1\\
\pm i & 0 &-1 & 0
\end{pmatrix}$$, see \cite{berndt} for more details.

Let us introduce also the complex basis vectors
$$Z = \begin{pmatrix} 0 & 0 & 1 & 0\\ 0 & 0 & 0 &1\\ -1 & 0 & 0 & 0\\ 0 & -1 & 0 & 0\end{pmatrix}, H'=-i \begin{pmatrix} 0 & 0 & 1 & 0\\ 0 & 0 & 0 &-1\\ -1 & 0 & 0 & 0\\ 0 & 1 & 0 & 0\end{pmatrix}$$
$$X=\frac{1}{2}\begin{pmatrix} 0 & 1 & 0 & -i\\ -1 & 0 & -i & 0\\ 0 & i & 0 & 1\\ i & 0 & -1 & 0\end{pmatrix}, \bar{X}=\frac{1}{2}\begin{pmatrix} 0 & 1 & 0 & i\\ -1 & 0 & i & 0\\ 0 & -i & 0 & 1\\ -i & 0 & -1 & 0\end{pmatrix}.$$
It is easy to check that
$$\mathfrak k_\mathbb C = \langle Z,H',X,\bar{X} \rangle_\mathbb C, \mathfrak h_\mathbb C = \langle Z,H'\rangle_\mathbb C, \mbox{ while }\mathfrak h = \langle H_1,H_2 \rangle_\mathbb R,$$ 
with commutation relations
$$\begin{matrix}
[Z,\mathfrak k_{\mathbb C}] &=& 0,\\
[H', X] &=& 2X,\\
[H',\bar{X}] &=& -2\bar{X},\\
[X,\bar{X}] &=& H'\end{matrix}$$

It means that the center of $\mathfrak k_\mathbb C$ is $C(\mathfrak k_\mathbb C) = \langle Z \rangle_\mathbb C$ and
$\mathfrak k_\mathbb C / C( \mathfrak k_\mathbb C) =\langle H',X,\bar{X}\rangle = \sl(2,\mathbb C)$.
$$\mathfrak p_\mathbb C = \langle Z,H',X_{(2,0)},X_{(0,2)},X_{(1,1)},X_{(1,-1)} \rangle_\mathbb C$$

There are the obvious relationship between the real root vectors and the complex basis of algebras as
$$\left\{\begin{array}{rcl}X_{\pm(2,0)} &=& \mp iH_1 + H_1 \pm 2iE_{2\mathbf e_1}\\
X_{\pm(1,1)} &=& \pm 2\bar{X} + 2E_{\mathbf e_1 -\mathbf e_2} \pm 2i E_{\mathbf e_1+\mathbf e_2}\\
X_{\pm(0,2)} &=& \pm H' + H_1 \pm 2i E_{2\mathbf e_2}
\end{array}\right.$$ 

Interchange the third and fourth basis vectors, the Cartan subgroup $H$ can be realized as 
$$H = \exp \mathfrak h = \left\{r(\theta_1)r(\theta_2) = \begin{pmatrix}
\cos\theta_1 & 0  & \sin\theta_1& 0\\
0 & \cos\theta_2  & 0& \sin\theta_2\\
-\sin\theta_1  & 0 &  \cos\theta_1 & 0\\
0 &  -\sin\theta_2&0 & \cos\theta_2 \end{pmatrix}\right\}.$$
The natural characters of this compact Cartan subgroup are realized as
$$r(\theta_1)r(\theta_2) \mapsto \exp \left(m_1\theta_1+m_2\theta_2 \right).$$
\begin{prop}
There are exactly two nontrivial endoscopic groups:
\begin{enumerate}
\item[a.] in the elliptic case $\mathbb S^1 \times \mathbb S^1 \times \{\pm 1\},$
\item[b.] in the parabolic case $\SL(2,\mathbb R) \times \{\pm 1\}.$
\end{enumerate}
\end{prop}
\textsc{Proof.}
Following the general rule, we do choose $\diag(\lambda_1 , \lambda_2, -\lambda_1,-\lambda_2)$  as a regular elliptic element of the Cartan subalgebra.
Then take the normalizers and take the connected component $Cent(\lambda, \mathfrak g)^0, \lambda = \lambda_1, \lambda_2, -\lambda_1,-\lambda_2)$.
There are 2 cases:
\begin{itemize}
\item[a.] two distinguished $\lambda_1 \ne \lambda_2$. In this case the connected component of identity in the normalizer is $\mathbb S^1 \times \mathbb S^1\times \{ \pm 1\}$, realized in form
$$\mathbb S^1 \times \mathbb S^1\times \{ \pm 1\}=\left\{\pm\begin{pmatrix}
\cos\theta_1 & 0  & \sin\theta_1& 0\\
0 & \cos\theta_2  & 0& \sin\theta_2\\
-\sin\theta_1  & 0 &  \cos\theta_1 & 0\\
0 &  -\sin\theta_2&0 & \cos\theta_2 \end{pmatrix}\right\}$$

\item[b.]  two identically equal $\lambda_1=\lambda_2$. In this case the connected component of the centralizer is $\SL(2,\mathbb R) \times \{\pm 1 \}$, realized in the form of 
$$\SL(2,\bbR) \times \{ \left. \pm 1\}=\left\{\pm\begin{pmatrix}
1 & 0  & 0& 0\\
0 & a  & b & 0\\
0 &c & d& 0\\
0 &  0&0 & 1 \end{pmatrix}\right| ad-bc =1\right\}$$
\end{itemize}
\hfill$\Box$

\subsection{Holomorphic Induction}
The discrete series representations of $\Sp(4,\mathbb R)$ is obtained by inducing from discrete series representations of two endoscopic groups.
Following the orbit method and the holomorphic induction, we do choose the integral functionals $\lambda$, take the corresponding orbits and then choose polarization and use the holomorphic induction.

As described above, the positive root system $\Delta^+ = \Delta^+_c \cup \Delta^+_n = \{(1,-1)\} \cup \{ (2,0),(1,1), (0,-2)\}, \rho = \frac{1}{2} \sum_{\alpha \in \Delta^+} \alpha = (2,-1)$. There are only two simple roots: compact one $(1,-1)$, and noncompact one $(0,-2)$. The corresponding coroots: $H_1=H_{1,-1}$ and $H_2=H_{0,-2}$ provides a basis of the Cartan subalgebra.

The discrete series representations of $\Sp(4,\bbR)$ were studied by many authors, I. Piateskii-Shapiro, R. Berndt, R. Berndt-W. Schmidt, etc .... and were decomposed into two series the representations $\sigma_k^+ = \ind_{\SL(2,\bbR) \ltimes \Heis(3,\bbR)}\pi^+_{m,k}$ $(m,k\in \bbN)$.
The other discrete series representations are obtained by cohomological induction as $\sigma_k^- = \cohind_{\SL(2,\bbR) \ltimes \Heis(3,\bbR)}\pi^-_{m,k}$ $(m,k\in \bbN)$ by R. Berndt \cite{berndt}.

We use the ordinary holomorphic induction to describe one part $\sigma_k^+$ of the discrete series.
The characters corresponding to this case are $\chi = (k,-k)$ and is reduced to the discrete series $\sigma^+_k$ of $\Sp(4,\bbR)$. The minimal $K$-type of $\sigma^+_k$ is $\tau_\Lambda$ and is $(k,-k)$,
$$\sigma^+_k = \Ind_{\SL(2,\bbR)\ltimes \Heis(3,\bbR)}^G\pi_{m,k}^+,$$ where
$$\pi^+_{m,k} = \pi^m_{SW} \otimes \pi^+_{k-1/2}$$ is the tensor product of the Shale-Weil representation $\pi^m_{SW}$ and the representation $ \pi^+_{k-1/2}$ of highest weight $\Lambda = (\lambda,\lambda'), \lambda \geq \lambda'$ and are integers.
Denote the highest weight by $\ell = \lambda + \lambda'$ and the lowest weight by $N = \lambda - \lambda'$, we have
$$\tau_\Lambda(e^{i\varphi}E_2) = e^{i\ell\varphi}E_{N+1},$$
$$\tau_\Lambda(t(\psi)) = \tau^\circ_N(t(\psi))= \diag(e^{-iN\psi},e^{iN\psi}), t(\psi) := \begin{pmatrix} e^{i\psi} & 0\\ 0 & e^{-i\psi} \end{pmatrix} $$
$\tau^\circ_\Lambda$ is the natural action of $g=t(\psi)$ on homogeneous polynomial of degree $N$,
$$\tau^\circ_N(g)P(\binom{u}{v}) = P(g^{-1}\binom{u}{v})$$ and therefore
$$\tau_\Lambda(\begin{pmatrix} e^{i\theta} & 0\\ 0& e^{-i\theta} \end{pmatrix}) = \diag\left(e^{i(\lambda'\theta + \lambda \theta)} , \dots, e^{i(\lambda'\theta + \lambda \theta)} \right)$$

\subsection{Cohomological Induction}
In this section we use cohomological to describe another part of the discrete series representations $\sigma_k^-$.
Let us consider the inclusion $\SU(1,1) \hookrightarrow \Sp(4,\bbC)$ in the natural way 
$$A= X+iY \mapsto \begin{pmatrix} X & -SY\\ SY & SXS  \end{pmatrix},$$ where $S = \begin{pmatrix} 1 & 0\\ 0 & -1 \end{pmatrix}$. The corresponding inclusion of Lie algebras is $$j: \su(1,1) = \langle U_1,U_2,U_3,U_4 \rangle\to \sp(4,\bbR),$$
$$j(U_1) = j(\begin{pmatrix} i & 0\\ 0 & 0 \end{pmatrix})= iZ = G-F = \begin{pmatrix} 0 & 0 & -1 & 0\\ 0 & 0 & 0 & 0\\ 1 & 0 & 0 & 0\\ 0 & 0 & 0 & 0   \end{pmatrix},$$
$$j(U_2) = j(\begin{pmatrix} 0 & 0\\ 0 & i \end{pmatrix})=iZ' = R-R' =\begin{pmatrix} 0 & 0 & 0 & 0\\ 0 & 0 & 0 & 1\\ 0 & 0 & 0 & 0\\ 0 & -1 & 0 & 0   \end{pmatrix},$$
$$j(U_3) = j(\begin{pmatrix} 0 & 1\\ 1 & 0 \end{pmatrix})=P_++P_- =\begin{pmatrix} 0 & 1 & 0 & 0\\ 1 & 0 & 0 & 0\\ 0 & 0 & 0 & -1\\ 0 & 0 & -1 & 0   \end{pmatrix},$$
$$j(U_4) = j(\begin{pmatrix} 0 & i\\ -i & 0 \end{pmatrix}) = -i(P_+-P_-) =\begin{pmatrix} 0 & 0 & 0 & 1\\ 0 & 0 & 1 & 0\\ 0 & 1 & 0 & 0\\ 1 & 0 & 0 & 0   \end{pmatrix}.$$
 The Lie algebra $\mathfrak l$ is exactly the centralizer of $$j(U_1) + j(U_2) = i(Z+Z') = iH'$$ and is $\theta$ stable,  thererefore $L$ is the stabilizer of this element in $\mathfrak g$.

Considqer the parabolic subgroup $Q \subset \Sp(4,\bbC)$, the Lie algebra of which is $$\mathfrak q = \mathfrak l + \mathfrak u = \langle Z,Z', P_\pm \rangle + \langle X_+, N_+, P_{0-} \rangle,$$
where the associate Levi subgroup $$L = \{ g\in G | \Ad(g) \mathfrak q \subset \mathfrak q \}$$, $$\mathfrak l = \Lie L \mbox{ and } \mathfrak l_\bbC = \mathfrak l \otimes \bbC \cong \su(1,1) = \sl(2,\bbC),$$
The Lie algebra $\mathfrak q$ is a polarization in the orbit method.

The unitary group $U(2)$ can be also included in the maximal compact subgroup $K$ of the group $\Sp(4,\bbR)$ by a map $j': \U(2) \to K \subset \Sp(4,\bbR)$ by 
$$A = X+iY \in \U(2) \mapsto \begin{pmatrix} X & -Y\\ Y & X\end{pmatrix}.$$ the corresponding map for the Lie algebras inclusion is
$$j': \mathfrak u(2) = \left\{\begin{pmatrix} \alpha & \beta + i\gamma\\ -\beta + i\gamma & i\delta \end{pmatrix} \right\}= \langle V_1,V_2,V_3,V_4 \rangle,$$
$$V_1 = U_1 + U_2= \begin{pmatrix} i & 0\\ 0 & i \end{pmatrix}, V_2 = U_1 - U_2 = \begin{pmatrix} i & 0\\ 0 & -i\end{pmatrix}$$
$$j'(U_1) =j(U_1) = j'(\begin{pmatrix} i & 0\\ 0 & 0 \end{pmatrix})= iZ = G-F = \begin{pmatrix} 0 & 0 & -1 & 0\\ 0 & 0 & 0 & 0\\ 1 & 0 & 0 & 0\\ 0 & 0 & 0 & 0   \end{pmatrix},$$
$$j'(U_2) =-j(U_2) = j'(\begin{pmatrix} 0 & 0\\ 0 & i \end{pmatrix})=-iZ' = R'-R =\begin{pmatrix} 0 & 0 & 0 & 0\\ 0 & 0 & 0 & -1\\ 0 & 0 & 0 & 0\\ 0 & 1 & 0 & 0   \end{pmatrix},$$
$$j'(V_3) = j'(\begin{pmatrix} 0 & 1\\ -1 & 0 \end{pmatrix})=N_++N_- =P'-P = =\begin{pmatrix} 0 & 1 & 0 & 0\\ -1 & 0 & 0 & 0\\ 0 & 0 & 0 & 1\\ 0 & 0 & 0 & -1   \end{pmatrix},$$
$$j'(V_4) = j'(\begin{pmatrix} 0 & i\\ i & 0 \end{pmatrix})=-i(N_+-N_-) =Q'-Q = =\begin{pmatrix} 0 &  & 0 & 0\\ -1 & 0 & -1 & 0\\ 0 & 1 & 0 & 0\\ 1 & 0 & 0 & 0   \end{pmatrix}.$$
Denote by $\bbT$ the compact Cartan subgroup, $\mathfrak t= \Lie \bbT$ ists Lie algebra, then it is easy to see that $$\mathfrak t = \mathfrak l \cap \mathfrak k = \mathfrak u(1,1) \cap \mathfrak u(2)$$ and
$$\mathfrak t_\bbC = \mathfrak t \otimes \bbC = \langle Z,Z'\rangle = \mathfrak h.$$
Because $T = L \cap K$ we have two fibrations
$$K/T \rightarrowtail Y = G/T \twoheadrightarrow X=\bbH=G/K \mbox{ and }  L/T \rightarrowtail Y \twoheadrightarrow D=G/L$$ with fibers $$K/T \cong U(2)/(U(1)\times U(1)) \cong SU(2)/U(1) \cong \bbP^1(\bbC)$$ and
$$L/T \cong U(1,1) /(U(1) \times U(1)) \cong \bbD^2 (\mbox{ open disc }), $$ respectively.

Following the cohomological induction, consider the Dolbeault cohomology of the complex of smooth (0,p) forms with values in the line bundle $\mathcal L_\chi$ with values in $\bbC_\chi$ of $\mathcal L_\chi$
$$A^p(D;\mathcal L_\chi) = \{C^\infty(G) \otimes \bbC_\chi \otimes \wedge^p\mathfrak u \}$$
The cohomology space $H^p(D;\mathcal L_\chi)$ is admissible representation of $G$ is is the maximal globalization of the Harish-Chandra module $\mathcal R^p_\mathfrak q(\bbC_{\chi -2\rho(\mathfrak u)}), $ where 
$\rho(\mathfrak u)= \frac{1}{2}(3,-3)$ is the half-sum of the positive roots from $\Delta(\mathfrak u) = \{(2,0),(1,-1),(0,-2) \}$. 

\begin{thm}
Let $s=\dim_\bbC K/(K\cap L) = 1$ be the dimension of a maximal compact subvariety of $D$ and $\chi$ is such a character of $T$, that
$$\label{cond} \langle \chi + \rho, \beta \rangle > 0, \forall \beta\in \Delta(\mathfrak u).$$ Then $$H^p(D; \mathcal L_\chi) =0, \forall p \ne s,$$
\end{thm}
Under the dominant condition \ref{cond}, Zierau described the Penrose transform
$$\mathcal P : H^s(D; \mathcal L_\chi) \to C^\infty(G/K,\mathcal E_{\chi'})$$ as an injection, where $\mathcal E_{\chi'}$ is the bundle over $G/K$ associated with the $K$-representation $E_{\chi'}$ of the fibers. Now following the Borel-Weil-Bott theorem for $H^s(K/(K\cap L); \check{\mathcal L}_{\check{\chi}})$, where $\check{\mathcal L}_{\check{\chi}}$ is the pull-back ove the holomorphic injection $K/(K\cap L) \hookrightarrow G/L$, we have the vanishing assertion. In our case, because the minimal $k$-type $\tau_\Lambda$ of $\sigma_k^-$ is $\lambda = (k-1,1-k)$, $k\geq 3$, $\lambda = \lambda -(2,-2) = (k-3, 3-k)$ and $\chi = \lambda + 2\rho(\mathfrak u) = \lambda +(3,-3) = (k,-k)$  we gave $$s=1, \chi=(k,-k), \chi' = (k-1,1-k).$$

\subsection{Hochschild-Serre spectral sequence}
Remark that because in general $\mathfrak p$ is not a subalgebra, we can modify it by taking subalgebra $\mathfrak h_+ = \mathbb C(Y+iX) \oplus \mathbb C(S-iZ/2)$:
$$\mathfrak e = \mathfrak p_+ \oplus \mathfrak k_\mathbb C = \mathfrak h_+ \oplus \mathfrak k_\mathbb C, \quad\mathfrak h_+.$$ Therefore, one has
$$\mathfrak e \cap \mathfrak b_1 = \mathfrak h_+, \quad \mathfrak e \cap \mathfrak b = \mathfrak h_+ \oplus \mathfrak m_\mathbb C.$$
We may construct a Hochschild-Serre spectral sequence for this filtration.

Consider a highest weight $\lambda + \alpha_{31}$ representation  $V^*_\lambda$ of $\mathfrak k_\mathbb C$, which is trivially on $\mathfrak p_+$ extended to a representation $\xi$ of $\mathfrak e = \mathfrak p_+ \oplus \mathfrak k_\mathbb C$. The action of $\mathfrak h_+$ in $V^{\lambda + \alpha_{31}}$ is $\xi + \frac{1}{2}\tr\ad_{\mathfrak b_1}$.
Denote by $\mathcal H_\pm$ the space of representations $T_\pm$ of $B$ $\Omega_\pm$ above and by $\mathcal H_\pm^\infty$ the subspaces of smooth vectors.
Because $\dim_\mathbb C(\mathfrak p_\mathbb C) = 2$, we have $\wedge^q(\mathfrak h_+)= 0$, for all $q \geq 3$.
It is natural to define the Hochschild-Serre cobound operators
$$(\delta_\pm)_{\lambda,q} : \wedge^q(\mathfrak h_+)^* \otimes V^{\lambda + \alpha_{31}} \otimes \mathcal H_\pm^\infty \to  \wedge^{q+1}(\mathfrak h_+)^* \otimes V^{\lambda + \alpha_{31}} \otimes \mathcal H_\pm^\infty$$ and by duality their formal adjoint operators $(\delta_\pm)_{\lambda,q}^*$.
The Hochschild-Serre spectral sequence is convergent
$$\bigoplus_{r+s=q} H^r(\mathfrak e_1; H^s(M;V^{\lambda + \alpha_{31}} \otimes \mathcal H_\pm^\infty)) \Longrightarrow H^q(B;\mathfrak b_1,V_\lambda)$$

\section{Trace Formula}
In this section we make precise the Arthur-Selberg trace formula for $\Sp(4,\bbR)$. 
\subsection{Trace formula}
Let us remind that $\Gamma \subset \Sp(4,\bbR)$ is a finitely generated Langlands type discrete subgroup with finite number of cusps.
Le $f\in C^\infty_c(\Sp(4,\bbR))$ be a smooth function of compact support. If $\varphi$ is a function from the representation space, the action of the induced representation $\inf_P^G\chi$ is the restriction of the right regular representation $R$ on the inducing space of induced representation.
$$\tr R(f)\varphi = \int_G (f(y)R(y)\varphi(x)dy) = \int_G f(y)\varphi(xy)dy$$
$$= \int_G f(x^{-1}y)\varphi(y)dy (\mbox{right invariance of Haar measure } dy) $$
$$= \int_{\Gamma\backslash G}\left(\sum_{\gamma\in \Gamma}f(x^{-1}\gamma y) 
\right)\varphi(y)dy$$

 Therefore, this action can be represented by an operator with kernel $K(x,y)$ of form
$$[R(f)\varphi](x) = \int_{\Gamma\backslash G} K_f(x,y)\varphi(y)dy,$$ where
$$K_f(x,y) = \sum_{\gamma\in \Gamma}f(x^{-1}\gamma y).$$ Because the function $f$ is of compact support, this sum convergent, and indeed is a finite sum, for any fixed $x$ and $y$ and is of class $L^2(\Gamma\backslash G \times\Gamma\backslash G )$. The operator is of trace class and it is well-known that
$$\tr R(f) = \int_{\Gamma\backslash G}K_f(x,x)dx.$$
As supposed, the discrete subgroup $\Gamma$ is finitely generated. Denote by $\{\Gamma\}$ the set of representatives of conjugacy classes. For any $\gamma\in\Gamma$ denote the centralizer of $\gamma\in \Omega \subset G$ by $\Omega_\gamma$, in particular, $G_\gamma \subset G$. Following the Fubini theorem for the doble integral, we can change the order of integration to have
$$\tr R(f) = \int_{\Gamma\backslash G} K_f(x,x)dx =\int_{\Gamma\backslash G} \sum_{\gamma\in \Gamma} f(x^{-1}\gamma x)dx$$
$$=\int_{\Gamma\backslash G} \sum_{\gamma\in \{\Gamma\}}\sum_{\delta\in \Gamma_\gamma\backslash \Gamma} f(x^{-1}\delta^{-1}\gamma\delta x)dx$$
$$=\sum_{\gamma\in \{\Gamma\}}\int_{\Gamma_\gamma\backslash G}  f(x^{-1}\gamma x)dx=\sum_{\gamma\in \{\Gamma\}}\int_{G_\gamma\backslash G}\int_{\Gamma_\gamma\backslash G_\gamma}  f(x^{-1}u^{-1}\gamma u x)dudx$$
$$= \sum_{\gamma\in \{\Gamma\}}\int_{G_\gamma\backslash G} \vol(\Gamma_\gamma \backslash G_\gamma)f(x^{-1}\gamma x)dx.$$
Therefore, in order to compute the trace formula, one needs to do:
\begin{itemize}
\item classfiy the conjugacy classes of all $\gamma$ in $\Gamma$: they are of type elliptic (different eigenvalues of the same sign),
 hyperbolic (nondegenerate, with eigenvalues of different sign),
 parabolic (denegerate)  
\item Compute the volume of form; it is the volume of the quotient of the stabilazer of the adjoint orbits. 
$\vol(\Gamma_\gamma \backslash G_\gamma)$
\item and compute the orbital integrals of form
$$\mathcal{O}(f) = \int_{G_\gamma\backslash G} f(x^{-1}\gamma x)d\dot x$$ 
The idea is to reduce these integrals to smaller endoscopic subgroups in order to the correponding integrals are ordinary or almost  ordinary. 
\end{itemize}

\subsection{Stable trace formula}
The Galois group $\Gal(\mathbb C/\mathbb R) = \mathbb Z_2$ of the complex field $\mathbb C$ is acting on the discrete series representation by character $\kappa(\sigma) = \pm 1$. Therefore the sum of characters can be rewrite as some sum over stable classes of characters.
$$\tr R(f) = \sum_{n=1}^\infty \sum_{\varepsilon = \pm 1} (\Theta_n^+(f) - \Theta_n^-(f)).$$

\section{Endoscopy}

\subsection{Orbital integrals}
\textsl{The simplest case} is the elliptic case when $\gamma = \diag(a_1,a_2,a_1^{-1}, a_2^{-1}) $. In this case, because of Iwasawa decomposition $x=mauk$, and the $K$-bivariance,  the orbital integral is
$$\mathcal O_\gamma(f) = \int_{G_\gamma\backslash G} f(x^{-1}\gamma x)dx = \int_U f(u^{-1}\gamma u)du$$ $$=\int_\mathbb R f(\begin{pmatrix} 1 & x & -y& z\\ 0 & 1& 0&y \\ 0 & 0 & 1 & x\\ 0 & 0 & 0 &  1\end{pmatrix}^{-1}\begin{pmatrix} a_1 & 0 & 0 & 0\\ 0 & a_2 & 0 & 0 \\ 0 & 0 & a_1^{-1} & 0\\ 0 & 0& 0 & a_2^{-1} \end{pmatrix} \begin{pmatrix} 1 & x & -y& z\\ 0 & 1& 0&y \\ 0 & 0 & 1 & x\\ 0 & 0 & 0 &  1\end{pmatrix})dsdxdydz$$ 
$$= |a_1-a_1^{-1}|^{-1}|a_2 -a_2^{-1}|^{-1} \mathcal O_\gamma(f).$$
The integral is abosolutely and uniformly convergent and therefore is smooth function of $a\in \bbR^*_+$. Therefore the function $$f^H(\gamma) = \Delta(\gamma)^{-1}\mathcal O_\gamma(f), \quad \Delta(\gamma) =|a_1-a_1^{-1}||a_2 -a_2^{-1}|$$ is a smooth function on the endoscopic group $H= (\bbR^*)^2$. 

\textsl{The second case} is the case where $\gamma = k_{\theta_1}k_{\theta_2} = \begin{pmatrix}\cos\theta_1 & 0 &\sin\theta_1 & 0\\ 0 & \cos\theta_2 &0 & \sin\theta_2 \\ -\sin\theta_1 &0 &  \cos\theta_1 & 0\\  0 & -\sin\theta_2 & 0 &\cos\theta_2 \end{pmatrix}$. We have again, $x= mauk$ and
$$\mathcal O_{k(\theta)}(f) = \int_{G_{k(\theta)}\backslash G} f(k^{-1}u^{-1}a^{-1}m^{-1}k(\theta)mauk)dmdudadk$$
$$= \int_{G_{k(\theta)}\backslash G} f(u^{-1}a^{-1}m^{-1}k(\theta)mau)dmduda$$
$$ = \int_{G_{k(\theta)}\backslash G} f(\begin{pmatrix} 1 & x  & -y & z \\ 0 & 1 & 0 & y \\  0 & 0 & 1 & x \\ 0 & 0 & 0 & 1  \end{pmatrix}^{-1}\begin{pmatrix} a_1^{-1} & 0 & 0 & 0\\ 0 & a_2^{-1}& 0 & 0\\ 0 & 0 & a_1 & 0\\  0 & 0 & 0 & a_2 \end{pmatrix} \begin{pmatrix}\cos\theta_1 & 0 &\sin\theta_1 & 0\\ 0 & \cos\theta_2 &0 & \sin\theta_2 \\ -\sin\theta_1 &0 &  \cos\theta_1 & 0\\  0 & -\sin\theta_2 & 0 &\cos\theta_2 \end{pmatrix}$$ $$\times \begin{pmatrix} a_1 & 0 & 0 & 0\\ 0 & a_2& 0 & 0\\ 0 & 0 & a_1^{-1} & 0\\  0 & 0 & 0 & a_2^{-1} \end{pmatrix}\begin{pmatrix} 1 & x  & -y & z \\ 0 & 1 & 0 & y \\  0 & 0 & 1 & x \\ 0 & 0 & 0 & 1  \end{pmatrix})duda dk(\theta_1)dk(\theta_2)$$
$$= \int_1^\infty \int_1^\infty f(\begin{pmatrix}\cos\theta_1 & 0 &t_1\sin\theta_1 & 0\\ 0 & \cos\theta_2 &0 & t_2\sin\theta_2 \\ -t_1^{-1}\sin\theta_1 &0 &  \cos\theta_1 & 0\\  0 & -t_2^{-1}\sin\theta_2 & 0 &\cos\theta_2 \end{pmatrix})\prod_{i=1}^2|t_i-t_i^{-1}|\frac{dt_i}{t_i}$$ 
$$= \int_0^{+\infty}\int_0^{+\infty}\sgn(t_1-1)\sgn(t_2-1) f(\begin{pmatrix}\cos\theta_1 & 0 &t_1\sin\theta_1 & 0\\ 0 & \cos\theta_2 &0 & t_2\sin\theta_2 \\ -t_1^{-1}\sin\theta_1 &0 &  \cos\theta_1 & 0\\  0 & -t_2^{-1}\sin\theta_2 & 0 &\cos\theta_2 \end{pmatrix})dt.$$ 
When $f$ is an element of the Hecke algebra, i.e. $f$ is of class $C^\infty_0(G)$ and is $K$-bivariant, the integral is converging absolutely and uniformly. Therefore the result is a function $F(\sin\theta)$. 
The function $f$ has compact support, then the integral is well convergent at $+\infty$. At the another point $0$, we develope the function $F$ into the Tayor-Lagrange of the first order with respect to $\lambda = \sin\theta \to 0$
$$F(\lambda) = A(\lambda) + \lambda B(\lambda),$$ where $A(\lambda) = F(0)$ and $B(\lambda)$ is the error-correction term $F'(\tau)$ at some intermediate value $\tau, 0 \leq \tau \leq t$.
Remark that 
$$ \begin{pmatrix}1 & 0 & 0 & 0\\ 0 & \sqrt{1-\lambda^2} & t\lambda & 0\\ 0 & -t^{-1}\lambda & \sqrt{1-\lambda^2}& 0 \\ 0 & 0 & 0 & 1 \end{pmatrix}$$ 
$$=\begin{pmatrix} 1 & 0 & 0 & \\ 0 & t^{1/2}&0& 0\\ 0 & 0 & t^{-1/2} & 0\\
 0 & 0 & 0 & 1 \end{pmatrix} \begin{pmatrix}1 & 0 & 0 & 0\\ 0 & \sqrt{1-\lambda^2} & t\lambda & 0\\ 0 & -t^{-1}\lambda & \sqrt{1-\lambda^2}& 0 \\ 0 & 0 & 0 & 1 \end{pmatrix}\begin{pmatrix} 1 & 0 & 0 & \\ 0 & t^{-1/2}&0& 0\\ 0 & 0 & t^{1/2} & 0\\
 0 & 0 & 0 & 1 \end{pmatrix} $$ we have
$$B= \frac{dF(\tau)}{d\lambda} =\frac{d}{d\lambda} \left. \int_0^{+\infty} \sgn(t-1)f(\begin{pmatrix}1 & 0 & 0 & 0\\ 0 & \sqrt{1-\lambda^2} & t\lambda & 0\\ 0 & -t^{-1}\lambda & \sqrt{1-\lambda^2}& 0 \\ 0 & 0 & 0 & 1 \end{pmatrix})dt \right|_{t=\tau}$$
$$=  \int_0^{+\infty} \sgn(t-1)g(\begin{pmatrix}1 & 0 & 0 & 0\\ 0 & \sqrt{1-\lambda^2} & t\lambda & 0\\ 0 & -t^{-1}\lambda & \sqrt{1-\lambda^2}& 0 \\ 0 & 0 & 0 & 1 \end{pmatrix}\frac{dt}{t}, $$
where $g\in C^\infty_c(N)$
and $g(\lambda) \cong O(-t^{-1}\lambda)^{-1}.$ $B$ is of logarithmic growth and $$B(\lambda) \cong \ln(|\lambda|^{-1})g(1)$$ up to constant term,  and therefore is contimuous. 
$$A = F(0) = |\lambda|^{-1} \int_0^\infty f(\begin{pmatrix}1 & 0 & 0 & 0\\ 0 &1 & \sgn(\lambda)u & 0\\ 0 & 0 & 1 & 0\\ 0 & 0 & 0 & 1\end{pmatrix}) du - 2f(I_3) + o(\lambda)$$
Hence the functions 
$$G(\lambda) = |\lambda|(F(\lambda) + F(\lambda)),$$
$$H(\lambda) = \lambda(F(\lambda) - F(-\lambda))$$ have the Fourier decomposition
$$G(\lambda) = \sum_{n=0}^N (a_n|\lambda|^{-1} + b_n)\lambda^{2n} + o(\lambda^{2N})$$
$$H(\lambda) = \sum_{n=0}^N h_n\lambda^{2n} + o(\lambda^{2N})$$
Summarizing the discussion, we have that in the case of $\gamma = k(\theta)$, there exists also a continuous function $f^H$ such that $$f^H(\gamma) = \Delta(\gamma) (\mathcal O_\gamma(f) - \mathcal O_{w\gamma}(f))=\Delta(k(\theta))\mathcal {SO}_\gamma(f), $$ where $\Delta(k(\theta)) = 4i\sin\theta_1\sin\theta_2$.

\begin{thm}
There is a natural function $\varepsilon : \Pi \to \pm 1$ such that in the Grothendieck group of discrete series representation ring, $$\sigma_G = \sum_{\pi\in \Pi} \varepsilon(\pi)\pi,$$ the map $\sigma \mapsto \sigma_G$ is dual to the map of geometric transfer, that for any $f$  on $G$, there is a unique $f^H$ on $H$
$$\tr \sigma_G(f) = \tr\sigma(f^H).$$
\end{thm}
\textsc{Proof.}
There is a natural bijection $\Pi_\mu \cong \mathfrak D(\mathbb R,H,G)$, we get a pairing
$$\langle .,.\rangle : \Pi_\mu \times \mathfrak k(\mathbb R,H,G)\to \mathbb C.$$ Therefore we have
$$\tr \Sigma_\nu(f^H) = \sum_{\pi\in \Pi_\Sigma} \langle s,\pi\rangle \tr \pi(f).$$
\hfill$\Box$

Suppose given a complete set of endoscopic groups $H = \mathbb S^1 \times \mathbb S^1 \times \{\pm 1\}$ or $\SL(2,\mathbb R) \times \{\pm 1\}$.  For each group, there is a natural inclusion
$$\eta: {}^LH \hookrightarrow {}^LG$$

Let $\varphi: DW_\bbR  \to {}^LG$ be the Langlands parameter, i.e. a homomorphism from the Weil-Deligne group $DW_\bbR = W_\bbR\ltimes \bbR^*_+$ the Langlands dual group,
$\bbS_\varphi$ be  the set of conjugacy classes of Langlands parameters modulo the connected component of identity map. For any $s\in 
\bbS_\varphi$,
$\check{H}_s = \Cent(s,\check G)^\circ$ the connected component of the centralizer of $s\in\bbS_\varphi$ we have $\check{H}_s$ is conjugate with $H$. 
Following the D. Shelstad pairing $$\langle s, \pi\rangle : \bbS_\varphi \times \Pi(\varphi) \to \bbC$$
$$\varepsilon(\pi) = c(s)\langle s,\pi\rangle.$$
Therefore, the relation
$$\sum_{\sigma\in \Sigma_s} \tr \sigma(f^H) = \sum_{\pi\in\Pi} \varepsilon(\pi) \tr \pi(f)$$
can be rewritten as
$$\widetilde{\Sigma}_s(f^H) = \sum_{s\in\Pi} \langle s,\pi\rangle \tr\pi(f)$$
and
$$\widetilde{\Sigma}_s(f^H) = c(s)^{-1}\sum_{\sigma\in\widetilde{\Sigma}_s} \tr \sigma(f^H) .$$
We arrive, finally to the result
\begin{thm}
$$\tr \pi(f) = \frac{1}{\#\bbS_\varphi}\sum_{s\in\bbS_\varphi} \langle s,\pi\rangle\widetilde{\Sigma}_s(\check{f}^H).$$
\end{thm}

\subsection{Stable orbital integral}

Let us remind that the \textit{orbital integral} is defined as
$$\mathcal{O}_\gamma(f) = \int_{G_\gamma\backslash G} f(x^{-1}\gamma x)d\dot x$$ 

The complex Weyl group is isomorphic to $\mathfrak S_3$ while the real Weyl group
is isomorphic to $\mathfrak S_2$ . The set of conjugacy classes inside a strongly regular
stable elliptic conjugacy class is in bijection with the pointed set
$\mathfrak S_3 /\mathfrak S_2$
 that can be viewed as a sub-pointed-set of the group
$\mathfrak E(\mathbb R, T, G) = (Z_2 )^2$
We shall denote by $\mathfrak K(\mathbb R, T, G)$ its Pontryagin dual.

Consider $\kappa \ne 1$ in $\mathfrak K(\mathbb R, T, G)$ such that $\kappa(H_{13}) = 1$. 
 Such a $\kappa$ is unique:
in fact one has necessarily
$\kappa(H_{12}) = \kappa(H_{13}) = -1$.

The endoscopic group $H$ one associates to $\kappa$ is isomorphic to
$S(U (1, 1) \times U (1))$
the positive root of $\mathfrak h$ in $H$ (for a compatible order) being $\alpha_{23} = \rho$

The endoscopic group $H$ can be embedded in $G$ as 
$$g(u,v,w) = \begin{pmatrix}
ua & iub & 0\\  -iuc & ud & 0\\  0 & 0 & v \end{pmatrix}, w = \begin{pmatrix} a & b\\ c & d \end{pmatrix}, ad-bc = 1 \mbox{ and } |u| = |v| = 1, uv = 1.$$
It will be useful to consider also the twofold cover
$H_1 = S(U (1) \times U (1)) \times \SL(2)$.

Let $f_\mu$ be a pseudo-coefficient for the discrete series representation $\pi_\mu$ then the \textit{$\kappa$-orbital integral} of a  regular element $\gamma$
in $T (R)$ is given by
$$\mathcal O^\kappa_\gamma(f_\mu) =\int_{G_\gamma\backslash G} \kappa(x) f_\mu(x^{-1}\gamma x) d\dot x = \sum_{\sign(w) =1} \kappa(w) \Theta^G_\mu(\gamma^{-1}_w)  = \sum_{\sign(w) =1} \kappa(w)\Theta_{w\mu}(\gamma^{-1}),$$ because there is a natural bijection between the left coset classes and the right coset classes.

\section{Poisson Summation Formula}
In the Langlands picture of the trace formula, the trace of the restriction of the regular representation on the cuspidal parabolic part is the coincidence of of the spectral side and the geometric side.
\begin{equation}
\sum_{\pi} m(\pi) \hat{f}(\pi) = \sum_{\gamma\in \Gamma\cap H} a^G_\gamma \hat{f}(\gamma) \end{equation}
Let us do this in  more details.

\subsection{Endoscopic transfer}
The transfer factor $\Delta(\gamma, \gamma_H )$ is given by
$$\Delta(\gamma, \gamma_H )=(−1)^{q(G)+q(H)} \chi_{G,H} (\gamma)\Delta_B (\gamma^{−1} ) . \Delta_{B_H}(\gamma_H^{-1})^{-1}$$
for some character $\chi_{G,H}$ defined as follows. Let $\xi$ be a character of the twofold covering $\mathfrak h_1$ of $\mathfrak h$, then
$\chi_{G,H} (\gamma^{−1} ) = e^{\gamma^{\rho - \rho_H +\xi}}$  defines a character of $H$, corresponding to $\mathfrak h$, because it is trivial on any fiber of the cover.

With such a choice we get when $\sign (w) = 1$
and $w\ne 1$, we have $\kappa(w) = -1$ and 
$$\Delta(\gamma^{-1},\gamma_H^{-1})\Theta_{w\mu}^G(\gamma) = -\frac{\gamma_H^{w\mu + \xi} - \gamma_H^{w_0w\mu + \xi}}{\gamma^{\rho_H}\Delta_{B_H}(\gamma_H}$$ therefore
$$\Delta(\gamma,\gamma_H)\Theta^G_{w\mu} (\gamma^{-1}) = \kappa(w)^{-1}\mathcal{SO}_\nu^H(\gamma_H^{-1}),$$ where $\nu = w\mu +\xi$ is running over the corresponding $L$-package of discrete series representations for the endoscopic group $H$.
Therefore we have the following formula
$$\Delta(\gamma,\gamma_H){\mathcal O}_\gamma^\kappa(f_\mu) = \sum_{\nu=w\mu + \xi\atop \sign(w) = 1} {\mathcal {SO}}_\nu^H(\gamma_H^{-1})$$ or
$$\Delta(\gamma,\gamma_H){\mathcal O}_\gamma^\kappa(f_\mu) = \sum_{\nu=w\mu + \xi\atop \sign(w) = 1} {\mathcal {SO}}_{\gamma_H}(g_\nu), $$ where $g_\nu$ is pseudo-coefficient for any one of the discrete series representation of the endoscopic subgroup $H$ in the $L$-package of $mu$.

For any $$f^H = \sum_{\nu=w\mu+\rho\atop \sign(w) = 1} a(w,\nu)g_\nu,\quad a(w_1,w_2\mu = \kappa(w_2) \kappa(w_2w_1)^{-1}$$ we have the formula
$$\tr \Sigma_\nu (f^H) = \sum_{w} a(w,\nu) \tr \pi_{w\mu}(f).$$

\subsection{Poisson summation and endoscopy}
\begin{thm}\cite{labesse}
There is a natural function $\varepsilon : \Pi \to \pm 1$ such that in the Grothendieck group of discrete series representation ring, $$\sigma_G = \sum_{\pi\in \Pi} \varepsilon(\pi)\pi,$$ the map $\sigma \mapsto \sigma_G$ is dual to the map of geometric transfer, that for any $f$  on $G$, there is a unique $f^H$ on $H$
$$\tr \sigma_G(f) = \tr\sigma(f^H).$$
\end{thm}
\textsc{Proof.}
There is a natural bijection $\Pi_\mu \cong \mathfrak D(\mathbb R,H,G)$, we get a pairing
$$\langle .,.\rangle : \Pi_\mu \times \mathfrak k(\mathbb R,H,G)\to \mathbb C.$$ Therefore we have
$$\tr \Sigma_\nu(f^H) = \sum_{\pi\in \Pi_\Sigma} \langle s,\pi\rangle \tr \pi(f).$$
\hfill$\Box$

Suppose given a complete set of endoscopic groups $H = \mathbb S^1 \times \mathbb S^1 \times \{\pm 1\}$ or $\SL(2,\mathbb R) \times \{\pm 1\}$.  For each group, there is a natural inclusion
$$\eta: {}^LH \hookrightarrow {}^LG$$

Let $\varphi: DW_\bbR  \to {}^LG$ be the Langlands parameter, i.e. a homomorphism from the Weil-Deligne group $DW_\bbR = W_\bbR\ltimes \bbR^*_+$ the Langlands dual group,
$\bbS_\varphi$ be  the set of conjugacy classes of Langlands parameters modulo the connected component of identity map. For any $s\in 
\bbS_\varphi$,
$\check{H}_s = \Cent(s,\check G)^\circ$ the connected component of the centralizer of $s\in\bbS_\varphi$ we have $\check{H}_s$ is conjugate with $H$. 
Following the D. Shelstad pairing $$\langle s, \pi\rangle : \bbS_\varphi \times \Pi(\varphi) \to \bbC$$
$$\varepsilon(\pi) = c(s)\langle s,\pi\rangle.$$
Therefore, the relation
$$\sum_{\sigma\in \Sigma_s} \tr \sigma(f^H) = \sum_{\pi\in\Pi} \varepsilon(\pi) \tr \pi(f)$$
can be rewritten as
$$\widetilde{\Sigma}_s(f^H) = \sum_{s\in\Pi} \langle s,\pi\rangle \tr\pi(f)$$
and
$$\widetilde{\Sigma}_s(f^H) = c(s)^{-1}\sum_{\sigma\in\widetilde{\Sigma}_s} \tr \sigma(f^H) .$$
We arrive, finally to the result
\begin{thm}\cite{labesse}
$$\tr \pi(f) = \frac{1}{\#\bbS_\varphi}\sum_{s\in\bbS_\varphi} \langle s,\pi\rangle \widetilde{\Sigma}_s(\check{f}).$$
\end{thm}

\begin{thm}
$$\tr R(f)_{L^2_{cusp}(\Gamma \backslash \Sp(4,\bbR)} = \sum_{\Pi_\mu}\sum_{\pi\in\Pi_\mu}m(\pi)\mathcal S\Theta_\pi(f) = \sum_{ \Pi_\mu} \Delta(\gamma,\gamma_H)\mathcal {SO}(f_\mu),$$ where
$$\mathcal S\Theta_\pi(f) = \sum_{\pi\in \Pi}\kappa(\pi) \Theta_\pi(f)$$ is the sum of Harish-Chandra characters of the discrete series running over the stable conjugacy classes of $\pi$ 
and $$\mathcal {SO}(f_\mu) = \sum_{\lambda\in \Pi_\mu}\kappa(\pi_\lambda) \mathcal O(f_\lambda)$$ is the sum of orbital integrals weighted by a character $\kappa : \Pi_\mu \to \{\pm1\}$. 
\end{thm}
\beginpf\;
The proof just is  a combination of the previous theorems.
\endpf

\end{document}